\input amstex
\documentstyle{amsppt}
\magnification=\magstep1
 \hsize 13cm \vsize 18.35cm \pageno=1
\loadbold \loadmsam
    \loadmsbm
    \UseAMSsymbols
\topmatter
\NoRunningHeads
\title $q$-Euler numbers and polynomials associated with
multiple $q$-zeta functions
\endtitle
\author
  T. Kim
\endauthor

\abstract The purpose of this paper is to present a systemic study
of some families of multiple $q$-Euler numbers and polynomials and
we construct multiple $q$-zeta function which interpolates
multiple $q$-Euler numbers at negative integer. This is a partial
answer of the open question in a previous publication( see J. Physics A: Math. Gen. 34(2001)7633-7638).
\endabstract
\thanks  2000 AMS Subject Classification: 11B68, 11S80
%\newline keywords and phrases
\newline  The present Research has been conducted by the research
Grant of Kwangwoon University in 2010
\endthanks
\endtopmatter

\document

{\bf\centerline {\S 1. Introduction/ Preliminaries}}

 \vskip 15pt
Let $p$ be a fixed odd prime number. Throughout this paper $\Bbb
Z_p$, $\Bbb Q_p$, $\Bbb C$ and $\Bbb C_p$ will, respectively, denote
the ring of $p$-adic rational integers, the field of $p$-adic
rational numbers, the complex number field and the completion of
algebraic closure of $\Bbb Q_p$. Let $v_p$ be the normalized
exponential valuation of $\Bbb C_p$ with
$|p|_p=p^{-v_p(p)}=\frac{1}{p}.$ When one talks of $q$-extension,
$q$ is variously considered  as and indeterminate, a complex number
$q\in \Bbb C$ or $p$-adic number $q\in \Bbb C_p$. If $q\in \Bbb C$,
one normally assumes $|q|<1$. If $q\in\Bbb C_p$, one normally
assumes $|1-q|_p<1.$ We use the notation
$$[x]_q=\frac{1-q^x}{1-q}, \text{ and }
[x]_{-q}=\frac{1-(-q)^x}{1+q}, \text{ (see [4, 5, 6, 7])}.$$ The
$q$-factorial is defined as $[n]_q!=[n]_q[n-1]_q\cdots[2]_q[1]_q$.
For a fixed  $d\in \Bbb N$ with $(p,d)=1$, $ d\equiv 1$ $ ( mod \
2)$ , we set
$$\aligned
&X =X_d =\varprojlim_N \Bbb Z/dp^N , X^\ast =\bigcup_{\Sb  0<a
<dp\\(a,p)=1
\endSb}a +d p\Bbb Z_p ,\\
& a +d p^N \Bbb Z_p = \{ x\in X | x \equiv a
\pmod{p^N}\},\endaligned$$ where $a\in\Bbb Z$ satisfies the
condition $0\leq a < d p^N$.
 The $q$-binomial formulae are known as
$$(b;q)_n=(1-b)(1-bq)\cdots(1-bq^{n-1})=\sum_{i=0}^n
{\binom{n}{i}}_q q^{\binom{i}{2}}(-1)^ib^i, $$ and
$$\frac{1}{(b;q)_n}=\frac{1}{(1-b)(1-bq)\cdots(1-bq^{n-1})}
=\sum_{i=0}^{\infty}{\binom{n+i-1}{i}}_q
b^i,$$ where
${\binom{n}{k}}_q=\frac{[n]_q!}{[n-k]_q![k]_q!}=\frac{[n]_q
[n-1]_q\cdots [n-k+1]_q}{[k]_q!}, \text{ (see [4, 8, 9])}.$

Recently, many authors have studied the $q$-extension in the
various area(see [4, 5, 6]). In this paper, we try to consider the
theory of $q$-integrals in the $p$-adic number field associated
with Euler numbers and polynomials closely related to fermionic
distribution. We say that $f$ is uniformly differentiable function
at a point $a\in\Bbb Z_p$, and write $f\in UD(\Bbb Z_p),$ if the
difference quotient $F_{f}(x, y)=\frac{f(x)-f(y)}{x-y}$ have a
limit $f^{\prime}(a)$ as $(x,y)\rightarrow (a,a).$ For $f \in
UD(\Bbb Z_p),$ the fermionic $p$-adic $q$-integral on $\Bbb Z_p$
is defined as
$$I_{q}(f)=\int_{\Bbb Z_p}f(x) d\mu_{q}(x)=\lim_{N\rightarrow
\infty}\frac{1+q}{1+q^{p^N}}\sum_{x=0}^{p^N-1}f(x)(-q)^x, \text{
(see [7, 8, 9])}. \tag1$$ Thus, we note that
$$\lim_{q\rightarrow 1}I_{q}(f)=I_{1}(f)=\int_{\Bbb Z_p} f(x)
d\mu_{1}(x). \tag2$$ For $n\in \Bbb N$, let $f_n(x)= f(x+n)$. Then
we have
$$I_{1}(f_n)=(-1)^nI_{1}(f)+2\sum_{l=0}^{n-1}(-1)^{n-1-l}f(l). \tag3$$

Using formula (3), we can readily derive  the Euler polynomials,
$E_n(x),$ namely,
$$\int_{\Bbb Z_p} e^{(x+y)t}d\mu_{1}(y)=\frac{2}{e^t
+1}e^{xt}=\sum_{n=0}^{\infty}E_n(x)\frac{t^n}{n!}, \text{ (see
[16-20])}.$$ In the special case $x=0$, the sequence
$E_n(0)=E_n$ are called the $n$-th Euler numbers. In one of an
impressive series of papers( see[1, 2, 3, 21, 23]), Barnes developed the
so-called multiple zeta and multiple gamma functions. Barnes'
multiple zeta function $\zeta_N (s, w| a_1, \cdots, a_N)$ depend
on the parameters $a_1, \cdots, a_N$ that will be assumed to be
positive. It is defined by the following series:
$$\zeta_{N}(s, w|a_1, \cdots, a_N)=\sum_{m_1, \cdots,
m_N=0}^{\infty}(w+m_1a_1+\cdots+m_Na_N)^{-s}\text{ for $\Re(s)>N,
\Re(w)>0 .$}\tag4$$ From (4), we can easily see that
$$\zeta_{M+1}(s, w+a_{M+1}|a_1, \cdots, a_{N+1})-\zeta_{M+1}(s, w|a_1,
\cdots, a_{N+1})=-\zeta_{M}(s, w|a_1, \cdots, a_N),$$
 and $\zeta_0(s, w)=w^{-s}$( see [11]). Barnes showed that
 $\zeta_N$ has a meromorphic continuation in $s$ (with simple
 poles only at $s=1, 2, \cdots, N$ and defined his multiple gamma
 function $\Gamma_N(w)$ in terms of the $s$-derivative at $s=0,$
 which may be recalled here as follows:   $\psi_n(w|a_1,\cdots, a_N)=
\partial_s\zeta_N(s,w|a_1, \cdots, a_N)|_{s=0}$( see[11]).
Barnes' multiple Bernoulli polynomials $B_n(x,r|a_1, \cdots, a_r)$
are defined by
$$\frac{t^r}{\prod_{j=1}^r(e^{a_jt}-1)}e^{xt}=\sum_{n=0}^{\infty}B_n(x,
r|a_1, \cdots, a_r)\frac{t^n}{n!}, \text{ ($|t|<2\pi$)}, \text{ (see [11])}.\tag5$$ By
(4) and (5), we see that
$$\zeta_N(-m, w|a_1, \cdots, a_N)=\frac{(-1)^N
m!}{(N+m)!}B_{N+m}(w,N|a_1, \cdots, a_N), $$ where $w>0$ and $m$
is a positive integer. By using the fermionic $p$-adic
$q$-integral on $\Bbb Z_p$, we consider the Barnes' type multiple
$q$-Euler polynomials and numbers in this paper. The main purpose
of this paper is to present a systemic study of some families of
Barnes' type multiple $q$-Euler polynomials and numbers. Finally,
we construct multiple $q$-zeta function which interpolates
multiple $q$-Euler numbers at negative integer. This is a partial
answer of the open question in [6, p.7637]

\vskip 10pt

{\bf\centerline {\S 2. Barnes type multiple $q$-Euler numbers and
polynomials }} \vskip 10pt

Let $x, w_1, w_2, \cdots, w_r$ be complex numbers with positive
real parts. In $\Bbb C$, the Barnes type multiple Euler numbers
and polynomials are defined by
$$\frac{2^r}{\prod_{j=1}^r(e^{w_jt}+1)}e^{xt}=\sum_{n=0}^{\infty}E_{n}^{(r)}(x|w_1,
\cdots, w_r)\frac{t^n}{n!}, \text{ for $|t|< max \{
\frac{\pi}{|w_i|}|i=1,\cdots,r \}$},\tag6$$ and $E_n^{(r)}(w_1,
\cdots, w_r)=E_n^{(r)}(0|w_1, \cdots, w_r)$(see [11, 12, 14]). In this
section, we assume that $q\in \Bbb C_p$ with $|1-q|_p<1$. We first
consider the $q$-extension of Euler polynomials as follows:
$$\sum_{n=0}^{\infty}E_{n,q}(x)\frac{t^n}{n!}=\int_{\Bbb
Z_p}e^{[x+y]_qt}d\mu_{1}(y)=2\sum_{m=0}^{\infty}(-1)^me^{[m+x]_qt},
\text{ (see [7, 8, 17]) } .\tag7$$ Thus, we have $E_{n,
q}(x)=\frac{2}{(1-q)^n}\sum_{l=0}^{n}\frac{\binom{n}{l}(-1)^lq^{lx}}{1+q^l}$(
see [7]). In the special case $x=0$, $E_{n,q}=E_{n,q}(0)$ are
called the $q$-Euler numbers. The $q$-Euler polynomials of order
$r\in\Bbb N$ are also defined by
$$\aligned
\sum_{n=0}^{\infty}E_{n,q}^{(r)}(x)\frac{t^n}{n!}&=\int_{\Bbb
Z_p}\cdots\int_{\Bbb
Z_p}e^{[x+x_1+\cdots+x_r]_qt}d\mu_{1}(x_1)\cdots d\mu_{1}(x_r)\\
&=2^r\sum_{m=0}^{\infty}\binom{m+r-1}{m}(-1)^m e^{[m+x]_qt},
\text{ (see [7, 8])}.
\endaligned\tag8$$
In the special case $x=0$, the sequence
$E_{n,q}^{(r)}(0)=E_{n,q}^{(r)}$ are refereed as the $q$-extension
of the Euler numbers of order $r$. Let $f\in\Bbb N$ with $f\equiv
1$ $(mod  \ 2)$. Then we have
$$\aligned
&E_{n,q}^{(r)}(x)=\int_{\Bbb Z_p}\cdots\int_{\Bbb
Z_p}[x+x_1+\cdots+x_r]_q^n d\mu_{1}(x_1)\cdots d\mu_{1}(x_r)\\
&=\frac{2^r}{(1-q)^n}\sum_{l=0}^{n}\binom{n}{l}(-1)^l q^{lx}
\sum_{a_1, \cdots, a_r=0}^{f-1}\sum_{m_1,\cdots,
m_r=0}^{\infty}q^{l(\sum_{i=1}^r(a_i+fm_i)}(-1)^{\sum_{i=1}^r
(a_i+fm_i)}\\
&=2^r\sum_{m_1, \cdots,
m_r=0}^{\infty}(-1)^{m_1+\cdots+m_r}[m_1+\cdots+m_r+x]_q^n.
\endaligned\tag9$$
By (8) and (9), we obtain the following theorem.

\proclaim{ Theorem 1} For $n\in \Bbb Z_+$, we have
$$\aligned
 E_{n,q}^{(r)}(x)&=2^r\sum_{m_1, \cdots,
m_r=0}^{\infty}(-1)^{m_1+\cdots+m_r}[m_1+\cdots+m_r+x]_q^n\\
&=2^r\sum_{m=0}^{\infty}\binom{m+r-1}{m}(-1)^m [m+x]_q^n.
\endaligned$$
\endproclaim
Let
$F_q^{(r)}(t,x)=\sum_{n=0}^{\infty}E_{n,q}^{(r)}(x)\frac{t^n}{n!}.$
Then we have
$$\aligned
F_q^{(r)}(t,x)&=2^r\sum_{m=0}^{\infty}\binom{m+r-1}{m}(-1)^me^{[m+x]_qt}\\
&=2^r\sum_{m_1, \cdots,
m_r=0}^{\infty}(-1)^{m_1+\cdots+m_r}e^{[m_1+\cdots+m_r+x]_qt}.
\endaligned\tag10$$
Let $\chi$ be the Dirichlet's character with conductor $f\in\Bbb
N$ with $f\equiv 1$ $( mod   \ 2)$. Then the generalized $q$-Euler
polynomials attached to $\chi$ are defined by
$$\sum_{n=0}^{\infty}E_{n,\chi,q}(x)\frac{t^n}{n!}=2\sum_{m=0}^{\infty}(-1)^m\chi(m)e^{[m+x]_qt}.
\tag11$$ Thus, we have
$$E_{n,\chi,q}(x)=\sum_{a=0}^{f-1}\chi(a)(-1)^a\int_{\Bbb
Z_p}[x+a+fy]_q^nd\mu_{1}(y)=[f]_q^n\sum_{a=0}^{f-1}\chi(a)(-1)^aE_{n,
q^f}(\frac{x+a}{f}).\tag12$$ In the special case $x=0$, the
sequence $E_{n,\chi,q}(0)=E_{n,\chi,q}$ are called the $n$-th
generalized $q$-Euler numbers attached to $\chi$. From (2) and
(3), we can easily derive the following equation.
$$E_{m,\chi,q}(nf)-(-1)^nE_{m,\chi,q}=2\sum_{l=0}^{nf-1}(-1)^{n-1-l}
\chi(l)[l]_q^m.$$ Let us consider higher-order generalized
$q$-Euler polynomials attached to $\chi$ as follows:
$$\int_X\cdots\int_X\left( \prod_{i=1}^r
\chi(x_i)\right)e^{[x_1+\cdots+x_r+x]_qt}d\mu_{1}(x_1)\cdots
d\mu_{1}(x_r)=\sum_{n=0}^{\infty}E_{n,\chi,q}^{(r)}(x)\frac{t^n}{n!},
\tag13$$ where $E_{n,\chi,q}^{(r)}(x)$ are called the $n$-th
generalized $q$-Euler polynomials of order $r$ attached to $\chi$.
By (13), we see that
$$\aligned
&E_{n,\chi,q}^{(r)}(x)=\frac{2^r}{(1-q)^n}\sum_{l=0}^n
\binom{n}{l}q^{lx}(-1)^l\sum_{a_1,\cdots,
a_r=0}^{f-1}\left(\prod_{j=1}^r\chi(a_j)\right)\frac{(-q^l)^{\sum_{i=1}^ra_i}}{(1+q^{lf})^r}\\
&=2^r\sum_{m=0}^{\infty}\binom{m+r-1}{m}(-1)^m\sum_{a_1,\cdots,
a_r=0}^{f-1}\left(\prod_{j=1}^r\chi(a_j)\right)(-1)^{\sum_{i=1}^r
a_i}[\sum_{j=1}^r a_j+x+mf]_q^n,
\endaligned\tag14$$
and
$$\sum_{n=0}^{\infty}E_{n,\chi,q}^{(r)}(x)\frac{t^n}{n!}
=2^r\sum_{m_1,\cdots,
m_r=0}^{\infty}(-1)^{\sum_{j=1}^rm_j}\left(\prod_{i=1}^r\chi(m_i)\right)
e^{[x+\sum_{j=1}^rm_j]_qt}. \tag15$$ In the special case $x=0$,
the sequence $E_{n,\chi,q}^{(r)}(0)=E_{n,\chi,q}^{(r)}$ are called
the $n$-th generalized $q$-Euler numbers of order $r$ attached to
$\chi$.

By (14) and (15), we obtain the following theorem.

\proclaim{ Theorem 2} Let $\chi$ be the Dirichlet's character with
conductor $f\in\Bbb N$ with $f\equiv 1$ $(mod  \ 2)$.
 For $n\in \Bbb Z_+$, $r\in\Bbb N$,  we have
$$\aligned
& E_{n,\chi,q}^{(r)}(x)=\frac{2^r}{(1-q)^n}\sum_{l=0}^n
\binom{n}{l}q^{lx}(-1)^l\sum_{a_1,\cdots,
a_r=0}^{f-1}\left(\prod_{j=1}^r\chi(a_j)\right)\frac{(-q^l)^{\sum_{i=1}^ra_i}}{(1+q^{lf})^r}\\
&=2^r\sum_{m=0}^{\infty}\binom{m+r-1}{m}(-1)^m\sum_{a_1,\cdots,
a_r=0}^{f-1}\left(\prod_{j=1}^r\chi(a_j)\right)(-1)^{\sum_{i=1}^r
a_i}[\sum_{j=1}^r a_j+x+mf]_q^n\\
&=2^r\sum_{m_1,\cdots,
m_r=0}^{\infty}(-1)^{m_1+\cdots+m_r}\left(\prod_{i=1}^r\chi(m_i)\right)
[x+m_1+\cdots+m_r]_q^n.
\endaligned$$
\endproclaim
  For $h\in\Bbb Z$ and $r\in\Bbb N$, we introduce the extended higher-order $q$-Euler polynomials
  as follows:
  $$E_{n,q}^{(h,r)}(x)=\int_{\Bbb Z_p}\cdots \int_{\Bbb Z_p}q^{\sum_{j=1}^r(h-j)x_j}[x+x_1+\cdots+x_r]_q^n
  d\mu_{1}(x_1)\cdots d\mu_1(x_r), \text{  (see [8]) }.\tag16$$
  From (16), we note that
  $$E_{n,q}^{(h,r)}(x)=2^r\sum_{m_1,\cdots, m_r=0}^{\infty}q^{(h-1)m_1+\cdots+(h-r)m_r}(-1)^{m_1+\cdots+m_r}
  [x+m_1+\cdots+m_r]_q^n.\tag17$$
  It is known in   [8] that
  $$ E_{n,q}^{(h,r)}(x)=\frac{2^r}{(1-q)^n}\sum_{l=0}^n \frac{\binom{n}{l}(-q^x)^l}{(-q^{h-r+l};q)_r}
    =2^r\sum_{m=0}^{\infty}{\binom{m+r-1}{m}}_q (-q^{h-r})^m[x+m]_q^n.\tag18$$                                                                                                Let $F_q^{(h,r)}(t,x)=\sum_{n=0}^{\infty}E_{n,q}^{(h,r)}(x)\frac{t^n}{n!}.$ Then we have
    $$\aligned
   F_q^{(h,r)}(t,x)&=2^r\sum_{m=0}^{\infty}{\binom{m+r-1}{m}}_qq^{(h-r)m}(-1)^me^{[m+x]_qt}\\
   &=2^r \sum_{m_1,\cdots, m_r=0}^{\infty}q^{\sum_{j=1}^r (h-j)m_j}(-1)^{\sum_{j=1}^m m_j}e^{[x+m_1+\cdots+m_r]_qt}.
    \endaligned\tag19$$
 Therefore, we obtain the following theorem.
  \proclaim{ Theorem 3} For $h, \in \Bbb Z$, $r\in\Bbb N$, and $x\in \Bbb Q^{+}$,
   we have
$$ E_{n,q}^{(h,r)}(x)=2^r\sum_{m_1, \cdots,
m_r=0}^{\infty}q^{(h-1)m_1+\cdots+(h-r)m_r}(-1)^{m_1+\cdots+m_r}[m_1+\cdots+m_r+x]_q^n.
$$
\endproclaim
For $f\in\Bbb N$ with $f\equiv 1$ $(mod   \ 2)$, it is easy to show that  the following distribution relation for $E_{n,q}^{(h,r)}(x)$.
 $$E_{n,q}^{(h,r)}(x)=[f]_q^n\sum_{a_1,\cdots, a_r=0}^{f-1}(-1)^{a_1+\cdots+a_r}q^{\sum_{j=1}^r(h-j)a_j}
  E_{n,q^f}(\frac{x+a_1+\cdots+a_r}{f}).$$
Let us consider Barnes' type multiple $q$-Euler polynomials.
For $w_1, \cdots, w_r \in\Bbb Z_p$, we define  the Barnes' type $q$-multiple Euler polynomials as follow:
$$E_{n,q}^{(r)}(x|w_1, \cdots, w_r)=\int_{\Bbb Z_p}\cdots \int_{\Bbb Z_p}[\sum_{j=1}^r w_jx_j+x]_q^n d\mu_{1}(x_1)\cdots d\mu_{1}(x_r).\tag20$$
From (20), we can easily derive the following equation.
$$ E_{n,q}^{(r)}(x|w_1, \cdots, w_r)=\frac{2^r}{(1-q)^n}\sum_{l=0}^n\frac{\binom{n}{l}(-q^x)^l}{(1+q^{lw_1})\cdots(1+q^{lw_r})}, \text{ (see [8] )}.
\tag21$$
Thus, we have
 $$E_{n,q}^{(r)}(x|w_1, \cdots, w_r)=\frac{2^r}{(1-q)^n}\sum_{l=0}^n\binom{n}{l}(-q^x)^l\sum_{a_1, \cdots, a_r=0}^{f-1}
 \frac{(-1)^{\sum_{i=1}^r a_i}q^{l\sum_{j=1}^r w_ja_j}}{(1+q^{lfw_1})\cdots(1+q^{lfw_r})},\tag22 $$
   where $f\in\Bbb N$ with $f\equiv 1   \  (mod    \   2)$.
   By (22), we see that
   $$ E_{n,q}^{(r)}(x|w_1, \cdots, w_r)=2^r\sum_{m_1,\cdots, m_r=0}^{\infty}(-1)^{m_1+\cdots+m_r}[x+w_1m_1+\cdots+w_rm_r]_q^n.\tag23$$
In the special case $x=0$, the sequence $E_{n,q}^{(r)}(w_1, \cdots, w_r)=E_{n,q}^{(r)}(0|w_1, \cdots, w_r)$ are called
the $n$-th Barnes' type multiple $q$-Euler numbers. Let $F_q^{(r)}(t,x|w_1, \cdots, w_r)=\sum_{n=0}^{\infty}E_{n,q}^{(r)}(x|w_1, \cdots, w_r)\frac{t^n}{n!}.$ Then we have
$$F_q^{(r)}(t,x|w_1, \cdots, w_r)=2^r\sum_{m_1, \cdots, m_r=0}^{\infty}(-1)^{m_1+\cdots+m_r}
e^{[x+w_1m_1+\cdots+w_rm_r]_qt}. \tag24$$
Therefore we obtain the following theorem.
   \proclaim{ Theorem 4} For $w_1, \cdots, w_r \in \Bbb Z_p$, $r\in\Bbb N$, and $x\in \Bbb Q^{+}$,
     we have
  $$\aligned
  E_{n,q}^{(r)}(x| w_1, \cdots, w_r)&=2^r\sum_{m_1, \cdots,
  m_r=0}^{\infty}(-1)^{m_1+\cdots+m_r}[x+m_1w_1+\cdots+m_rw_r]_q^n \\
 & =\frac{2^r}{(1-q)^n}\sum_{l=0}^n\frac{\binom{n}{l}(-q^x)^l}{(1+q^{lw_1})\cdots(1+q^{lw_r})}.
 \endaligned$$
  \endproclaim
  For $w_1, \cdots, w_r \in \Bbb Z_p$, $a_1, \cdots, a_r \in \Bbb Z$, we consider another $q$-extension
  of Barnes' type multiple $q$-Euler polynomials as follows:
 $$E_{n,q}^{(r)}(x|w_1,\cdots, w_r;a_1,\cdots, a_r)
 =\int_{\Bbb Z_p}\cdots \int_{\Bbb Z_p}[x+\sum_{j=1}^rw_jx_j]_q^nq^{\sum_{i=1}^r a_ix_i}\left(\prod_{i=1}^r  d\mu_{1}(x_i)\right).\tag25$$
 Thus, we have
   $$E_{n,q}^{(r)}(x|w_1,\cdots,w_r;a_1,\cdots, a_r)
   =\frac{2^r}{(1-q)^n}\sum_{l=0}^n  \frac{\binom{n}{l}(-1)^lq^{lx}}{(1+q^{lw_1+a_1})\cdots(1+q^{lw_r+a_r}) }.\tag26$$
 From (25) and (26), we can derive the following equation.
 $$E_{n,q}^{(r)}(x|w_1,\cdots, w_r;a_1, \cdots, a_r)
 =2^r\sum_{m_1,\cdots, m_r=0}(-1)^{\sum_{j=1}^rm_j}q^{\sum_{i=1}^r a_im_i}[x+\sum_{j=1}^rw_jx_j]_q^n .\tag27$$
 Let $F_q^{(r)}(t,x|w_1,\cdots, w_r;a_1,\cdots, a_r)
 =\sum_{n=0}^{\infty}E_{n,q}^{(r)}(x|w_1,\cdots, w_r;a_1, \cdots, a_r)\frac{t^n}{n!}.$ Then, we have
 $$\aligned
 &F_q^{(r)}(t,x|w_1,\cdots, w_r;a_1,\cdots, a_r)\\
 &=2^r\sum_{m_1,\cdots, m_r=0}^{\infty}
 (-1)^{m_1+\cdots+m_r}q^{ a_1m_1+\cdots+a_rm_r}e^{[x+w_1m_1 +\cdots+ w_rm_r
 ]_qt}.\endaligned\tag28$$

   \proclaim{ Theorem 5} For  $r\in\Bbb N$, $w_1, \cdots, w_r\in\Bbb Z_p$, and $a_1, \cdots, a_r \in \Bbb Z$,
   we have
 $$ E_{n,q}^{(r)}(x|w_1, \cdots, w_r; a_1, \cdots a_r)=2^r\sum_{m_1, \cdots,
 m_r=0}^{\infty}(-1)^{\sum_{j=1}^rm_j}q^{\sum_{i=1}^r a_im_i}[x+\sum_{j=1}^r w_jm_j]_q^n.$$
 \endproclaim
 Let $\chi$ be a Dirichlet's character with conductor $f\in \Bbb N$ with $f\equiv 1$ $(mod  \  2)$.
 Now we consider the generalized Barnes' type $q$-multiple Euler polynomials attached to $\chi$ as follows:
 $$\aligned
 &E_{n,\chi,q}^{(r)}(x|w_1, \cdots, w_r;a_1, \cdots, a_r)\\
& =\int_X\cdots \int_X [x+w_1x_1+\cdots+w_rx_r]_q^n\left(\prod_{j=1}^r\chi(x_j)\right)q^{a_1x_1+\cdots +a_rx_r}d\mu_{1}(x_1)\cdots d\mu_{1}(x_r)
.\endaligned$$
Thus, we have
$$\aligned
    &E_{n, \chi,q}^{(r)}(x|w_1,\cdots, w_r;a_1, \cdots, a_r) \\
    &=\frac{2^r}{(1-q)^n}\sum_{b_1, \cdots, b_r=0}^{f-1}\left(\prod_{i=1}^r\chi(b_i)\right)(-1)^{\sum_{j=1}^rb_j}
   q^{\sum_{i=1}^r(lw_i+a_i)b_i}
   \frac{\sum_{l=0}^n \binom{n}{l}(-1)^l q^{lx}}{\prod_{j=1}^r(1+q^{(lw_j+a_j)f})}. \endaligned\tag 29$$
   From, (29), we note that
    $$\aligned
    &E_{n, \chi,q}^{(r)}(x|w_1,\cdots, w_r;a_1, \cdots, a_r) \\
    &=2^r\sum_{m_1, \cdots, m_r=0}^{\infty}\left(\prod_{j=1}^r\chi(m_i)\right)(-1)^{m_1+\cdots+m_r}
    q^{a_1m_1+\cdots+a_rm_r}[x+\sum_{j=1}^rw_jm_j]_q^n    . \endaligned$$
Therefore we obtain the following theorem.
  \proclaim{ Theorem 6} For  $r\in\Bbb N$, $w_1, \cdots, w_r\in\Bbb Z_p$, and $a_1, \cdots, a_r \in \Bbb Z$,
   we have
 $$\aligned
    &E_{n, \chi,q}^{(r)}(x|w_1,\cdots, w_r;a_1, \cdots, a_r) \\
    &=2^r\sum_{m_1, \cdots, m_r=0}^{\infty}\left(\prod_{j=1}^r\chi(m_i)\right)(-1)^{m_1+\cdots+m_r}
    q^{a_1m_1+\cdots+a_rm_r}[x+\sum_{j=1}^rw_jm_j]_q^n    . \endaligned$$
 \endproclaim
 Let $F_{q, \chi}^{(r)}(t,x|w_1,\cdots, w_r;a_1,\cdots, a_r)
  =\sum_{n=0}^{\infty}E_{n,\chi, q}^{(r)}(x|w_1,\cdots, w_r;a_1, \cdots, a_r)\frac{t^n}{n!}.$

  By Theorem 6, we see that
  $$\aligned
   &F_{q, \chi}^{(r)}(t,x|w_1,\cdots, w_r;a_1,\cdots, a_r) \\
   &= 2^r\sum_{m_1, \cdots, m_r=0}^{\infty}\left(\prod_{j=1}^r\chi(m_i)\right)(-1)^{m_1+\cdots+m_r}
      q^{a_1m_1+\cdots+a_rm_r}e^{[x+\sum_{j=1}^rw_jm_j]_qt}.    \endaligned\tag30$$

\vskip 10pt
   {\bf\centerline {\S 3. Barnes type multiple $q$-zeta functions }} \vskip 10pt

   In this section, we assume that $q\in\Bbb C$ with $|q|<1$ and the parameters $w_1,
  \cdots, w_r$ are positive.  From (28), we consider the Barnes' type multiple $q$-Euler polynomials in $\Bbb C$ as follows:
  $$\aligned
   &F_q^{(r)}(t, x|w_1, \cdots, w_r;a_1, \cdots, a_r)\\
   &= 2^r\sum_{m_1,\cdots, m_r=0}^{\infty}(-1)^{m_1+\cdots+m_r}q^{a_1m_1+\cdots+a_rm_r}  e^{[x+w_1m_1+\cdots+w_rm_r]_qt}\\
   &=\sum_{n=0}^{\infty} E_{n,q}^{(r)}(x|w_1, \cdots, w_r;a_1,\cdots, a_r)\frac{t^n}{n!},
   \text{ for $|t|<\max_{1\leq i \leq r} \{ \frac{\pi}{|w_i|}\}$ }.
 \endaligned\tag31$$
     For $s, x \in \Bbb C$ with $\Re(x)>0,$ $a_1, \cdots, a_r \in\Bbb C$,
  we can derive the following Eq.(32) from the  Mellin transformation of $F_q^{(r)}(t,x|w_1, \cdots, w_r; a_1, \cdots, a_r)$.
  $$\aligned
  &\frac{1}{\Gamma(s)}\int_{0}^{\infty}t^{s-1}F_q^{(r)}(t,x|w_1, \cdots, w_r;a_1,\cdots, a_r)dt\\
 & = \sum_{m_1, \cdots, m_r=0}^{\infty}\frac{(-1)^{m_1+\cdots+m_r}q^{ m_1a_1+\cdots+m_ra_r}}{[x+ w_1m_1+\cdots+w_rm_r]_q^s}.
 \endaligned \tag32$$
  For $s, x \in \Bbb C$ with $\Re(x)>0$, $a_1, \cdots, a_r \in \Bbb C$,  we define Barnes' type multiple $q$-zeta function as follows:
  $$\zeta_{q,r}(s,x|w_1, \cdots, w_r;a_1, \cdots, a_r)
  =\sum_{m_1, \cdots, m_r=0}^{\infty}\frac{(-1)^{m_1+\cdots+m_r}q^{m_1a_1+\cdots+m_ra_r}}{[x+w_1m_1+\cdots+w_rm_r]_q^s}
 .\tag33$$
 Note that
 $\zeta_{q, r}(s,x|w_1,\cdots, w_r)$  is meromorphic function in whole complex $s$-plane. By using the Mellin transformation and the Cauchy
 residue theorem, we obtain the following theorem which is a part of answer of open question in [6, p.7637 ] .

  \proclaim{ Theorem 7} For $x \in \Bbb C$ with $\Re(x)>0$, $n\in\Bbb Z_{+}$, we have
   $$\zeta_{q,r}(-n, x|w_1, \cdots, w_r; a_1, \cdots, a_r)
   =(-1)^{n+1}E_{n,q}^{(r)}(x|w_1, \cdots, w_r;a_1, \cdots, a_r).$$
    \endproclaim
    Let $\chi$ be a Dirichlet's character with conductor $f\in \Bbb N$ with $f\equiv 1$ $(mod  \  2)$.
 From (30), we can define the generalized Barnes' type multiple $q$-Euler polynomials attached to $\chi$ in $\Bbb C$ as follows:
 $$\aligned
   &F_{q, \chi}^{(r)}(t,x|w_1,\cdots, w_r;a_1,\cdots, a_r) \\
   &= 2^r\sum_{m_1, \cdots, m_r=0}^{\infty}\left(\prod_{j=1}^r\chi(m_i)\right)(-1)^{m_1+\cdots+m_r}
      q^{a_1m_1+\cdots+a_rm_r}e^{[x+\sum_{j=1}^rw_jm_j]_qt}\\
      &= \sum_{n=0}^{\infty} E_{n,\chi,q}^{(r)}(x|w_1, \cdots, w_r;a_1,\cdots, a_r)\frac{t^n}{n!}.    \endaligned\tag34$$

       From (34) and  Mellin transformation of $F_{q, \chi}^{(r)}(t,x|w_1, \cdots, w_r; a_1, \cdots, a_r)$,
        we can easily derive the following equation (35) .
  $$\aligned
  &\frac{1}{\Gamma(s)}\int_{0}^{\infty}t^{s-1}F_{q, \chi}^{(r)}(t,x|w_1, \cdots, w_r;a_1,\cdots, a_r)dt\\
 & = \sum_{m_1, \cdots, m_r=0}^{\infty}\frac{\left( \prod_{j=1}^r \chi(m_i)\right)(-1)^{m_1+\cdots+m_r}q^{ m_1a_1+\cdots+m_ra_r}}{[x+ w_1m_1+\cdots+w_rm_r]_q^s}.
 \endaligned \tag35$$
  For $s, x \in \Bbb C$ with $\Re(x)>0$,  we also define Barnes' type multiple $q$-$l$ function as follows:
  $$\aligned
  &l_{q,\chi}^{(r)}(s,x|w_1, \cdots, w_r;a_1, \cdots, a_r)\\
 & =\sum_{m_1, \cdots, m_r=0}^{\infty}\frac{\left(\prod_{j=1}^r\chi(m_j)\right)(-1)^{m_1+\cdots+m_r}q^{m_1a_1+\cdots+m_ra_r}}{[x+w_1m_1+\cdots+w_rm_r]_q^s}
 .\endaligned\tag36$$
    Note that
 $l_{q, \chi}^{(r)}(s,x|w_1,\cdots, w_r)$  is meromorphic function in whole complex $s$-plane.
 By using (34), (35), (36), and the Cauchy
 residue theorem, we obtain the following theorem.

  \proclaim{ Theorem 8} For $x, s \in \Bbb C$ with $\Re(x)>0$, $n\in\Bbb Z_{+}$, we have
   $$l_{q,\chi}^{(r)}(-n, x|w_1, \cdots, w_r; a_1, \cdots, a_r)
   =(-1)^{n+1}E_{n,\chi, q}^{(r)}(x|w_1, \cdots, w_r;a_1, \cdots, a_r).$$
    \endproclaim
 We note that Theorem 8 is $r$-multiplication of Dirichlet's type $q$-$l$-series.  Theorem 8 seems to be  interesting and worthwhile
 for doing study in the area of  multiple $p$-adic $l$-function or mathematical physics related to Knot theory and $\zeta$-function (see [4-20]).

\vskip 20pt

\quad Taekyun Kim

\quad Division of General Education-Mathematics,

\quad Kwangwoon University,

\quad Seoul 139-701, S. Korea \quad e-mail:\text{
tkkim$\@$kw.ac.kr}

\enddocument